\documentclass[12]{article}

\begin{document}

\title{Rigorous Covariant Path Integrals} 

\author{ Alexander Dynin \\ 
Ohio State University, Columbus, Ohio 43210, USA \\ E-mail 
dynin@math.ohio-state.edu}

\maketitle

\begin{abstract}
 Our rigorous path integral is extended to quantum 
evolution on metric-affine manifolds.
\end{abstract}

\section{Path Integrals on Euclidean Spaces.} Consider the evolution 
operator
$U[\psi (t',q)] = \psi (t'',q),\ t''\geq t'$,\ of the quantum 
evolution equation on $\mathcal L^2 (\mathbf R^d)$: $$ 
\frac{\hbar}{i}\partial\psi (t,q)/\partial t 
+f(t,q,\frac{\hbar}{i}\partial/\partial q)\psi (t,q)=0. $$ Here the 
operator $f(t,q,\frac{\hbar}{i} \partial /\partial q)$ is the 
standard quantization of a classical time-dependent quasi-Hamiltonian
$f(t,q,p)$ on the phase space $\mathbf R^{2d}$. 

The formal Hamiltonian functional integral $$
\int \prod_{t''\geq t\geq t'} \frac{ dq(t) dp(t)}{(2\pi \hbar)^{d}} 
(t) \exp \frac{i}{\hbar}
\int_{t'}^{t''} [p(t)\dot q(t)-f(t,q(t),p(t))]dt, $$ presumably 
represents the standard symbol $\mathcal U^{(st)}=\langle 
q|U|p\rangle \langle p|q\rangle$ of the evolution operator $U$.
In practice, such formal integral is just a convenient notation 
subjected to all standard algorithms of the integral calculus, 
including the
stationary phase approximations and analytic continuation to the 
imaginary time.
Application of other ordering rules of canonical quantization, such 
as Weyl, Wick and more general $\Omega$-rules~\cite{aga}, leads to so 
called $\Omega$-symbols~\cite{dyn}.

Actually, $\Omega=\Omega (q,p)$ is a real analytic function on the 
phase space
$\mathbf R^{2d}$. Suppose, to start with, that $\Omega (q,p)$ has no 
zeroes. Then, by definition, the $\Omega$-symbol of the operator 
$f(t,q,\frac{\hbar}{i} \partial /\partial q)$ is $$ f^{\Omega}(q,p)
= (1/\Omega) (-i\hbar \partial _{q},-i\hbar \partial 
_{p})f^{(we)}(q,p) $$
on $\mathbf R^{2d}$. Here $f^{(we)}(q,p)$ is the Weyl symbol of that 
operator.

For example, $\Omega$ is $\exp(\frac{i}{2} qp)$ for the standard 
quantization, and for the Wick quantization $\Omega$ is 
$\exp(-\frac{1}{2}(p^{2}+q^{2}))$.

In general, $\Omega (q,p)$ has real
zeroes. Then the $\Omega$-symbol exists only $\hbar 
$-asymptotically~\cite{dyn}.

An $\Omega$-symbol $f(q,p)$ is called a \emph{quasi-Hamiltonian} if, 
uniformly
in \mbox{$ t'\leq t \leq t''$,} it is

(i) \emph{quasi-polynomial} of order $m>0$: for any multi-index 
$\alpha$ there is a constant $C$ such that $ |\partial 
^{\alpha}f(q,p)|\leq C(1+\sqrt {q^{2}+p^{2}})^{m-|\alpha |}. $

(ii) \emph{quasi-dissipative}: $\Re (if) >\delta$, where $\delta$ a 
constant.

(iii) \emph{hypo-elliptic}~\cite{shub}.

(iv) \emph{continuous in} $t$ along with all its derivatives 
$\partial _{q,p}$.

Simple examples are the Schr\"{o}dinger Hamiltonians with arbitrary 
quasi-polynomial
scalar and vector potentials.

The path integral construction is based on the following 
\emph{Convergence Theorem}~\cite{dyn}:

Consider the ordered products $U_{\mathcal P}=\prod_{n=1}^{N}U_{n}$, 
where $U_{n}$ are the operators with the Weyl symbols 
$[1+\frac{i}{\hbar}f(t_{n},q,p)(t_{n}-t_{n-1})]^{-1}$, associated with
partitions $\mathcal P: t'=t_{0}<t_{1}<\ldots< t_{n-1}<t_{N}=t''$ of 
the time-interval $[t',t'']$.
Then the strong (i.e. in expectation) operator limit 
$U=\lim_{|\mathcal P|\rightarrow \infty}U_{\mathcal P}$ exists at 
least for small $\hbar >0$.

The existence of the operator limit entails the representation of the 
(not
quasi-polynomial) symbols $\mathcal U^{\Omega}$ as the limits of the 
quasi-polynomial
$\Omega$-symbols of $\mathcal U_{\mathcal P}^{\Omega}$. According to 
the $\Omega$-symbols product rules~\cite{dyn} the latter are multiple
distributional integral transforms of the products of the $\mathcal 
U_{n}^{\Omega}$.

For example, the standard symbol $\mathcal U_{\mathcal P}^{(st)}(q,p)$
is the limit of the multiple distributional integrals $$
\int_{\mathbf R^{2dN}} \prod _{n=1}^{N-1}\frac{dq_{n}dp_{n}}{(2\pi 
\hbar )^{dN}} \prod _{n=1}^{N} \mathcal U_{n}^{(st)}(q_{n},p_{n-1}) 
\exp\left [\frac{i}{\hbar}\sum_{n=1}^{N}p_{n-1}\cdot 
(q_{n}-q_{n-1})\right ]. $$
with the boundary conditions $q_{0}= q_{N}=q,\quad p_{0} =p_{N}=p$. 

If $\Omega (q,p)$ has zeroes, then
the functional integral is $\hbar$-asymptotic. 

\emph{By our definition, the path integrals are such limits. The 
integral calculus algorithms for such path integrals are justified by 
such approximations.}

For example, the spectacular reappearance of the classical 
Hamiltonian equations
from the semi-classical Feynman functional integral takes the 
following form:
the principal terms in the $\hbar $-expansions of the partial 
products symbols are
the backward Euler approximations of the corresponding classical 
Hamilton equations.

\section{Covariant Path Integrals on Manifolds.} 

The proof of the Convergence Theorem~\cite{dyn} is based on the theory
of the finite difference solutions of the Cauchy problem with 
necessary estimates,
provided by the calculus of the $\Omega$-symbols. Consequently, an 
extension of the functional integral to quantum evolution on manifolds
requires a covariant extension of the $\Omega$-calculus.

Let $Q$ be a finite-dimensional configuration manifold 
with a Riemannian metric $g^{ij}$, and $\nabla$ a (possibly 
non-symmetric)
metric connection on $Q$. Such metric-affine manifolds appear, e.g., 
in the Kleinert's functional integral for the hydrogen atom via a 
non-holonomic change of variables~\cite{klein}.

We assume that the connection has \emph{bounded geometry}, i.e. the 
curvature and torsion tensors and their covariant derivative of any 
order are uniformly
bounded relative to the given metric $g^{ij}$.

Let $\mathcal S(Q)$ and $\mathcal S'(Q)$ be the Schwartz spaces of 
the test functions
and temperate distributions (alias ket- and bra-vectors), defined as 
for the Euclidean $Q=\mathbf R^{d}$~\cite{shub}: $\mathcal S(Q)$ is 
the space of the functions with all covariant derivatives falling off 
faster than any negative power of the geodesic distance from a fixed 
observation point $o$, and $\mathcal S'(Q)$ is the anti-dual space of 
the anti-linear functionals on it, which are continuous relative to 
the natural Frechet topology of the $\mathcal S(Q)$. The anti-duality 
is assumed in the form
$
\langle \Psi | \psi \rangle = \int _{Q} dq\sqrt{g(q)}\Psi (q)^{*}\psi 
(q).
$

The $\mathcal S'(Q)$ is supplied with the strong topology of uniform 
convergence on bounded subsets
of $\mathcal S(Q)$.
Actually the Schwartz spaces do not depend on the connection of 
bounded geometry and the observation point. 

A continuous linear operator $A:\mathcal S(Q) \rightarrow \mathcal 
S'(Q)$ is called a \emph Schwartz operator. Any differential operator 
on $Q$ with coefficients from $\mathcal S'(Q)$ is such an operator. 
The celebrated \emph{Schwartz Kernel Theorem} states that 

(i) there is one to one correspondence between the Schwartz operators 
$A$ and their Schwartz kernels $\mathcal A(q'',q')\in S'(Q''\times 
Q')$ such that
$
\langle A\psi (q'') |\psi (q') \rangle =\langle \mathcal A(q'',q') | 
\psi (q'')\psi (q')\rangle,
$
and vice versa.

(ii) The strong
convergence of Schwartz operators is equivalent to the strong 
convergence of their Schwartz kernels.

Let $\mathcal N$ be a closed \emph{symmetric convex normal 
neighborhood} of the
diagonal of the square $Q\times Q$ with the induced metric and 
connection, so that for all $(q',q'') \in \mathcal N$ 

(i) $(q'',q') \in \mathcal N$,

(ii) there is a unique $\nabla $-geodesic $\gamma _{q'q''}$ from $q'$ 
to $q''$,

(iii) the inverse exponential
mapping $\exp_{q'}^{-1}(q'')$
defines
the normal coordinates on the normal neighborhood $ 
N_{q'}=\{q'':(q',q'')\in \mathcal N\}$.

We require our Schwartz operators to be $\mathcal N$-proper: their 
kernels should vanish outside of $\mathcal N$. Such are, for example, 
the differential operators with quasi-polynomial coefficients. 

The \emph{covariant Weyl symbol} of $A$ is defined (cf.~\cite{ful}) 

$
\mathcal A^{(we)}(q,p) = \int_{T_{q}} dv\sqrt{g(q)} e^{-ip\cdot v} 
\mathcal A (\exp _{q}(v/2), \exp _{q}(-v/2)),\quad p\in T^{*}_{q}. $

The covariant $\Omega$-symbol is defined via the Weyl symbol as $$ 
\mathcal A^{\Omega}(q,p) = (1/\Omega )(\frac{\hbar }{i}) (\nabla 
_{q}^{sym},\partial _{p}) \mathcal A^{(we)} (q,p), $$ where the order 
of the differentiation is such that the symmetric covariant $\nabla 
_{q}^{sym}$ are always prior to the partial $\partial _{p}$.

The definition of the covariant $\Omega$-symbol is compatible with 
that of Safarov ,~\cite{saf}
for his $s$-symbols, $0\leq s \leq 1$. In particular, $s=1$ 
corresponds to the standard symbol $
\mathcal A^{(st)}(q,p)= \int_{T_{q}} dv\sqrt{g(q)}\cdot e^{-ip\cdot 
v}\cdot \mathcal A (q, \exp _{q}v).
$
The quasi-dissipativity, hypo-ellipticity and $t$-continuity for 
quasi-Hamiltonian symbols are covariantly defined as in the Euclidean 
case.

\emph{The convergence theorem is valid on the the metric-affine 
manifolds $Q$ as well. The covariant
$\Omega$-symbols of an evolution operator are the limits of covariant 
multiple integrals over $(T^{*}Q)^{N}$.} In particular, the standard 
symbol $\mathcal U^{(st)}(q,p)$ is the strong limit
of the standard symbols $\mathcal U_{\mathcal P}^{(st)}(q,p)$, which 
are $$
\int _{(T^{*}Q)^{N}}d\lambda _{N}
\left [\prod _{n=1}^{N} \mathcal
U_{n}^{(st)}(q_{n},\exp_{q_{n}}^{-1}(p_{n-1})) \right ]\cdot \exp 
\left (\frac{i}{\hbar}\sum_{n=1}^{N}p_{n-1}\cdot \exp_{q_{n-1}}^{-1} 
(q_{n})\right ).
$$
Here
$
d\lambda _{N}=\prod _{n=1}^{N-1}\frac{dq_{n}dp_{n}}{(2\pi \hbar 
)^{dN}},\quad p_{n}\in T_{q_{n}}^{*}, \quad q_{0}= q_{N}=q,\quad 
p_{0} =p_{N}=p. $

For the functional integrals of evolution operators in metric vector 
bundles $E$ over $Q$ with connections $\nabla ^{E}$ of bounded 
geometry (such as spinor bundles over compact manifolds), the 
covariant kernels $\mathcal A(q'',q')$
and their $\Omega$-symbols are linear transformations of the fibers 
$E_{q'}\rightarrow E_{q''}$ (as for the Dirac operators). 
Accordingly, all corresponding formulas involve the \emph{parallel 
transport $\tau ^{E}_{q'q}$ along} $\gamma _{q'q''}$ in $E$. 

\section{Conclusion.} Our rigorous path integral is extended to quantum 
evolution on metric-affine manifolds

\end{document}